\documentclass[12pt]{amsart}
\usepackage{epsfig,amssymb,euscript}
\calclayout
\newcommand\R{\mathbb R}
\renewcommand\P{\mathcal P}
\newcommand\Q{\mathcal Q}
\renewcommand\S{\mathcal S}
\newcommand\T{\EuScript T}
\newcommand\Map{\mathsf S}
\newcommand\Aff{\operatorname{\mathsf{Aff}}}
\newcommand\Bil{\operatorname{\mathsf{Bil}}}
\newcommand\x{\times}
\newcommand\seminorm[1]{\lvert#1\rvert}
\newcommand\norm[1]{\lVert#1\rVert}

\theoremstyle{plain}
\newtheorem{thm}{Theorem}

\theoremstyle{remark}
\newtheorem*{remark}{Remark}
\newtheorem*{remarks}{Remarks}

\begin{document}

\title{Approximation by Quadrilateral Finite Elements}
\author{Douglas N. Arnold}
\address{Department of Mathematics, Penn State University, University Park,
PA 16802}
\email{dna@psu.edu}
\urladdr{http://www.math.psu.edu/dna/}
\author{Daniele Boffi}
\address{Dipartimento di Matematica, Universit\`e Pavia, 27100 Pavia, Italy}
\email{boffi@dimat.unipv.it}
\urladdr{http://dimat.unipv.it/\char'176boffi/}
\author{Richard S. Falk}
\address{Department of Mathematics, Rutgers University, Piscataway, NJ
08854}
\email{falk@math.rutgers.edu}
\urladdr{http://www.math.rutgers.edu/\char'176falk/}
\date{February 25, 2000}
\begin{abstract}
We consider the approximation properties of finite element spaces on
quadrilateral meshes.  The finite element spaces are constructed
starting with a given finite dimensional space of functions on a square
reference element, which is then transformed to a space of functions on
each convex quadrilateral element via a bilinear isomorphism of the
square onto the element.  It is known that for affine isomorphisms, a
necessary and sufficient condition for approximation of order $r+1$ in $L^2$
and order $r$ in $H^1$ is that the given space of functions on the reference
element contain all polynomial functions of total degree at most $r$.  In
the case of bilinear isomorphisms, it is known that the same estimates
hold if the function space contains all polynomial functions of
separate degree $r$.  We show, by means of a counterexample, that this
latter condition is also necessary. As applications we demonstrate
degradation of the convergence order on quadrilateral meshes as
compared to rectangular meshes for serendipity finite elements and for
various mixed and nonconforming finite elements.
\end{abstract}

\maketitle

\section{Introduction}
Finite element spaces are often constructed starting with a finite
dimensional space $\hat V$ of shape functions given on a reference
element $\hat K$ and a class $\Map$ of isomorphic mappings of the
reference element. If $F\in\Map$ we obtain a space of functions
$V_F(K)$ on the image element $K=F(\hat K)$ as the compositions of
functions in $\hat V$ with $F^{-1}$.  Then, given a partition $\T$
of a domain $\Omega$ into images of $\hat K$ under mappings in
$\Map$, we obtain a finite element space as a
subspace\footnote{The subspace is typically determined by some
interelement continuity conditions.  The imposition of such conditions
through the association of local degrees of freedom is an important
part of the construction of finite element spaces, but, not being
directly relevant to the present work, will not be discussed.} of the
space $V^{\T}$ of all functions on $\Omega$ which restrict to
an element of $V_F(K)$ on each $K\in\T$.

For example, if the reference element $\hat K$ is the unit
triangle, and the reference space $\hat V$ is the space $\P_r(\hat K)$
of polynomials of degree at most $r$ on $\hat K$, and the mapping class
$\Map$ is the space $\Aff(\hat K)$ of affine isomorphisms of $\hat
K$ into $\R^2$, then $V^{\T}$ is the familiar space of all
piecewise polynomials of degree at most $r$ on an arbitrary triangular
mesh $\T$.  When $\Map = \Aff(\hat K)$, as in this case,
we speak of \emph{affine finite elements}.

If the reference element $\hat K$ is the unit square, then it is often
useful to take $\Map$ equal to a larger space than $\Aff(\hat K)$,
namely the space $\Bil(\hat K)$ of all bilinear isomorphisms of $\hat
K$ into $\R^2$.  Indeed, if we allow only affine images of the unit
square, then we obtain only parallelograms, and we are quite limited as
to the domains that we can mesh (e.g., it is not possible to mesh
a triangle with parallelograms).  On the other hand, with bilinear
images of the square we obtain arbitrary convex quadrilaterals, which
can be used to mesh arbitrary polygons.

The above framework is also well suited to studying the approximation
properties of finite element spaces.  See, e.g., \cite{ciarlet-raviart}
and \cite{ciarlet}. A fundamental result holds in the case of affine
finite elements: $\Map=\Aff(\hat K)$. Under the assumption that the
reference space $\hat V\supseteq \P_r(\hat K)$,  the following result
is well known: if $\T_1$, $\T_2$, \ldots is any shape-regular sequence
of triangulations of a domain $\Omega$ and $u$ is any smooth function
on $\Omega$,  then the $L^2$ error in the best approximation of $u$ by
functions in $V^{\T_n}$ is $O(h^{r+1})$ and the piecewise $H^1$ error
is $O(h^r)$, where $h=h(\T_n)$ is the maximum element diameter. It is
also true, even if less well-known, that the condition that $\hat V\supseteq
\P_r(\hat K)$ is necessary if these estimates are to hold.

The above result does not restrict the choice of reference element
$\hat K$, so it applies to rectangular and parallelogram meshes by
taking $\hat K$ to be the unit square.  But it does not apply to
general quadrilateral meshes, since to obtain them we must choose
$\Map=\Bil(\hat K)$, and the result only applies to affine finite
elements. In this case there is a standard result analogous to the
positive result in the previous paragraph, \cite{ciarlet-raviart},
\cite{ciarlet}, \cite[Section I.A.2]{girault-raviart}. Namely, if $\hat
V\supseteq \Q_r(\hat K)$, then for any shape-regular sequence of
quadrilateral partitions of a domain $\Omega$ and any smooth function
$u$ on $\Omega$, we again obtain that  the error in the best
approximation of $u$ by functions in $V^{\T_n}$ is $O(h^{r+1})$ in
$L^2$ and  $O(h^r)$ in (piecewise) $H^1$.  It turns out, as we shall
show in this paper, that the hypothesis that $\hat V\supseteq \Q_r(\hat
K)$ is strictly necessary for these estimates to hold.  In particular,
if $\hat V\supseteq  \P_r(\hat K)$ but $\hat V\nsupseteq \Q_r(\hat K)$,
then the rate of approximation achieved on general shape-regular
quadrilateral meshes will be strictly lower than is obtained using
meshes of rectangles or parallelograms.

More precisely, we shall exhibit in Section~3 a domain $\Omega$ and a
sequence, $\T_1$, $\T_2$, \ldots of quadrilateral meshes of it,
and prove that whenever $V(\hat K)\nsupseteq \Q_r(\hat K)$, then there
is a function $u$ on $\Omega$ such that
\begin{equation*}
\inf_{v\in V^{\T_n}}\|u-v\|_{L^2(\Omega)}\ne o(h^r),
\end{equation*}
 (and so, {\it a fortiori}, is $\ne O(h^{r+1})$).
A similar result holds for $H^1$ approximation.  The counterexample is
far from pathological. Indeed, the domain $\Omega$ is as simple as
possible, namely a square; the mesh sequence $\T_n$ is simple and as
shape-regular as possible in that all elements at all mesh levels are
similar to a single fixed trapezoid; and the function $u$ is as smooth
as possible, namely a polynomial of degree $r$.

The use of a reference space which contains $\P_r(\hat K)$ but not
$\Q_r(\hat K)$ is not unusual, but the degradation of convergence order
that this implies on general quadrilateral meshes in comparison to
rectangular (or parallelogram) meshes is not widely appreciated.
It has been observed in special cases, often as a result of numerical
experiments, cf.~\cite[Section~8.7]{zienkiewicz}.

We finish this introduction by considering some examples. Henceforth we
shall always use $\hat K$ to denote the unit square.
First, consider finite elements with the simple
polynomial spaces as shape functions: $\hat V=\P_r(\hat K)$. These of
course yield $O(h^{r+1})$ approximation in $L^2$ for rectangular
meshes.  However, since $\P_r(\hat K)\supseteq \Q_{\lfloor
r/2\rfloor}(\hat K)$ but $\P_r(\hat K)\nsupseteq \Q_{\lfloor
r/2\rfloor+1}(\hat K)$, on general quadrilateral meshes they only
afford $O(h^{\lfloor r/2\rfloor+1})$ approximation.

A similar situation holds for serendipity finite element spaces, which have
been popular in engineering computation for thirty years.  These spaces
are constructed using as reference shape functions the space $\S_r(\hat
K)$ which is the span of $\P_r(\hat K)$ together with the two monomials
$\hat x^r\hat y$ and $\hat y\hat x^r$.  (The purpose of the additional
two functions is to allow local degrees of freedom which can be used to
ensure interelement continuity.)  For $r=1$, $\S_1(\hat K)=\Q_1(\hat K)$, but
for $r>1$ the situation is similar to that for $\P_r(\hat K)$, namely
$\S_r(\hat K)\supseteq \Q_{\lfloor r/2\rfloor}(\hat K)$ but $\S_r(\hat
K)\nsupseteq \Q_{\lfloor r/2\rfloor+1}(\hat K)$.  So, again, the
asymptotic accuracy achieved for general quadrilateral meshes is only
about half that achieved for rectangular meshes: $O(h^{\lfloor
r/2\rfloor+1})$ in $L^2$ and $O(h^{\lfloor r/2\rfloor})$ in $H^1$.  In
Section~4 we illustrate this with a numerical example.

While the serendipity elements are commonly used for solving second
order differential equations, the pure polynomial spaces $\P_r$ can
only be used on quadrilaterals when interelement continuity is not
required.  This is the case in several mixed methods.  For example, a
popular element choice to solve the stationary Stokes equations is
bilinearly mapped piecewise continuous $\Q_2$ elements for the two
components of velocity, and discontinuous piecewise linear elements for
the pressure. Typically the pressure space is taken to be functions
which belong to $\P_1(K)$ on each element $K$.  This is known to be a
stable mixed method and gives second order convergence in $H^1$ for the
velocity and $L^2$ for the pressure.  If one were to define the
pressure space instead by using the construction discussed above,
namely by composing linear functions on reference square with bilinear
mappings, then the approximation properties of mapped $\P_1$ discussed
above would imply that method could be at most first order accurate, at
least for the pressures.  Hence, although the use of mapped $\P_1$ as
an alternative to unmapped $\P_1$ pressure elements is sometimes
proposed \cite{shopov-iordanov}, it is probably not advisable.

Another place where mapped $\P_r$ spaces arise is for approximating the
scalar variable in mixed finite element methods for second order
elliptic equations.  Although the scalar variable is discontinuous, in
order to prove stability it is generally necessary to define the space
for approximating it by composition with the mapping to the
reference element (while the space for the vector variable is defined
by a contravariant mapping associated with the mapping to the reference
element). In the case of the Raviart--Thomas rectangular elements, the
scalar space on the reference square is $\Q_r(\hat K)$, which maintains
full $O(h^{r+1})$ approximation  properties under bilinear mappings. By
contrast, the scalar space used with the Brezzi-Douglas-Marini and the
Brezzi-Douglas-Fortin-Marini spaces is $\P_r(\hat K)$.  This necessarily
results in a loss of approximation order when mapped to quadrilaterals
by bilinear mappings.

Another type of element which shares this difficulty is the simplest
nonconforming quadrilateral element, which generalizes to
quadrilaterals the well-known piecewise linear non-conforming element
on triangles, with degrees of freedom at the midpoints of edges. On the
square, a bilinear function is not well-defined by giving its value at
the midpoint of edges (or its average on edges), because these
quantities do not comprise a unisolvent set of degrees of freedom (the
function $(\hat x-1/2)(\hat y-1/2)$ vanishes at the four midpoints of the edges
of the unit square).  Hence, various definitions of nonconforming
elements on rectangles replace the basis function $\hat x\hat y$ by some other
function such as $\hat x^2-\hat y^2$.  Consequently, the reference
space contains $\P_1(\hat K)$, but does not contain $\Q_1(\hat K)$, and
so there is a degradation of convergence on quadrilateral meshes.  This
is discussed and analyzed in the context of the Stokes problem in
\cite{rannacher-turek}.

As a final application, we remark that many of the finite element
methods proposed for the Reissner-Mindlin plate problem are based
on mixed methods for the Stokes equations and/or for second order
elliptic problems.  As a result, many of them suffer from the same sort
of degradation of convergence on quadrilateral meshes.  An analysis of
a variety of these elements will appear in forthcoming work by the
present authors.

In Section~3, we prove our main result, the necessity of the condition
that the reference space contain $\Q_r(\hat K)$ in order to obtain
$O(h^{r+1})$ approximation on quadrilateral meshes. The proof relies on
an analogous result for affine approximation on rectangular meshes,
where the space $\P_r(\hat K)$ enters rather than $\Q_r(\hat K)$. While
this is a special case of known results, for the convenience of the
reader we include an  elementary proof in Section~2. In the final
section we illustrate the results with numerical computations.

\section{Approximation theory of rectangular elements}
In this section we prove some results concerning approximation
by rectangular elements which will be needed in the next section
where the main results are proved.  The results in this
section are essentially known, and many are true in far greater
generality than stated here.

If $K$ is any square with edges parallel to the axes, then $K=F_K(\hat K)$
where $F_K(\hat x):=x_K + h_K\hat x$ with $x_K\in\R^2$  and $h_K>0$ the
side length. For any function $u\in L^2(K)$, we define
$\hat u_K= u\circ F_K\in L^2(\hat K)$, i.e.,
$\hat u_K(\hat x)=u(x_K+h_K\hat x)$.
Given a subspace $\hat S$ of $L^2(\hat K)$ we define the associated
subspace on an arbitrary square $K$ by
\begin{equation*}
S(K)=\{\, u:K\to\R\,|\, \hat u_K\in \hat S\,\}.
\end{equation*}
Finally, let $\Omega$ denote the unit cube ($\Omega$ and $\hat K$ both
denote the unit cube, but we use the notation $\Omega$ when we think of
it as a fixed domain, while we use $\hat K$ when we think of it as a
reference element).  For $n=1,2,\ldots$, let $\T_h$ be the uniform mesh
of $\Omega$ into $n^d$ subcubes when $h=1/n$, and define
\begin{equation*}
S_h = \{\, u:\Omega\to\R\,|\, u|_K\in S(K)\text{ for all $K\in\T_h$}\,\}.
\end{equation*}
In this definition, when we write $u|_K\in S(K)$ we mean only that $u|_K$
agrees with a function in $S_K$ almost everywhere, and so do not impose
any interelement continuity.

The following theorem gives a set of equivalent conditions for optimal
order approximation of a smooth function $u$ by elements of $S_h$.

\begin{thm}\label{thm:poly}
Let $\hat S$ be a finite dimensional subspace of $L^2(\hat K)$, 
$r$ a non-negative integer.   The following conditions are equivalent:
\begin{enumerate}
\item There is a constant $C$ such that
$\displaystyle\inf_{v\in S_h}\norm{u-v}_{L^2(\Omega)} \leq
 Ch^{r+1}\seminorm{u}_{r+1}$
for all $u\in H^{r+1}(\Omega)$.
\item $\displaystyle\inf_{v\in S_h}\norm{u-v}_{L^2(\Omega)} = o(h^r)$
for all $u\in \P_r(\Omega)$.
\item $\P_r(\hat K)\subset\hat S$.
\end{enumerate}
\end{thm}

\begin{proof}
The first condition implies that $\inf_{v\in S_h}\norm{u-v}_{L^2(\Omega)}=0$
for $u\in \P_r(\Omega)$, and so
implies the second condition.  The fact that the third condition
implies the first is a well-known consequence of the Bramble--Hilbert
lemma.  So we need only show that the second condition implies the
third.

The proof is by induction on $r$.  First consider the case $r=0$.  We have
\begin{equation}
\label{eq:sum}
\inf_{v\in S_h}\norm{u-v}_{L^2(\Omega)}^2
= \sum_{K\in\T_h}\inf_{v_K\in S(K)}\norm{u-v_K}_{L^2(K)}^2
= h^2\sum_{K\in\T_h}\inf_{w\in\hat S}\norm{\hat u_K-w}_{L^2(\hat K)}^2,
\end{equation}
where we have made the change of variable $w=\hat v_K$ in the last
step.

In particular, for $u\equiv1$ on $\Omega$, $\hat u_K\equiv1$ on $\hat K$ for all $K$, so the quantity
\begin{equation*}
c:=\inf_{w\in\hat S}\norm{\hat u_K-w}_{L^2(\hat K)}^2
\end{equation*}
is independent of $K$.  Thus
\begin{equation*}
\inf_{v\in S_h}\norm{u-v}_{L^2(\Omega)}^2
=h^2\sum_{K\in\T_h}c= c
\end{equation*}
The hypothesis that this quantity is $o(1)$ implies that $c=0$, i.e.,
that the constant function belongs to $\hat S$.

Now we consider the case $r>0$.  We again apply \eqref{eq:sum}, this
time for $u$ an arbitrary homogeneous polynomial of degree $r$.  Then
\begin{equation}\label{eq:hom}
\hat u_K(\hat x)=u(x_K+h\hat x)= u(h\hat x)+p(\hat x)
 = h^r u(\hat x)+p(\hat x),
\end{equation}
where $p\in\P_{r-1}(\hat K)$.  Substituting in \eqref{eq:sum},
and invoking the inductive hypothesis that $\hat S\supseteq\P_{r-1}(\hat K)$, we get that
\begin{equation*}
\inf_{v\in S_h}\norm{u-v}_{L^2(\Omega)}^2
=h^{2+2r}\sum_{K\in\T_h}
\inf_{w\in\hat S}\norm{u-w}_{L^2(\hat K)}^2
=h^{2r}\inf_{w\in\hat S}\norm{u-w}_{L^2(\hat K)}^2.
\end{equation*}
Again the last infimum is independent of $K$ so we immediately deduce that
$u$ belongs to $\hat S$.  Thus $\hat S$ contains all homogeneous
polynomials of degree $r$ and all polynomials of degree less than $r$
(by induction), so it indeed contains all polynomials of degree at most
$r$.
\end{proof}

A similar theorem holds for gradient approximation.  Since the
finite elements are not necessarily continuous we write
$\nabla_h$ for the gradient operator applied piecewise
on each element.

\begin{thm}\label{thm:poly1}
Let $\hat S$ be a finite dimensional subspace of $L^2(\hat K)$, 
$r$ a non-negative integer.   The following conditions are equivalent:
\begin{enumerate}
\item There is a constant $C$ such that
$\displaystyle\inf_{v\in S_h}\norm{\nabla_h(u-v)}_{L^2(\Omega)} \leq
Ch^r\seminorm{u}_{r+1}$
for all $u\in H^{r+1}(\Omega)$.
\item $\displaystyle\inf_{v\in S_h}\norm{\nabla_h(u-v)}_{L^2(\Omega)} =
o(h^{r-1})$
for all $u\in \P_r(\Omega)$.
\item $\P_r(\hat K)\subset\P_0(\hat K)+\hat S$.
\end{enumerate}
\end{thm}

\begin{proof}
Again, we need only prove that the second condition implies the third.
In analogy to \eqref{eq:sum}, we have
\begin{equation}
\label{eq:sum1}
\begin{aligned}
\inf_{v\in S_h}\sum_{K\in\T_h}\norm{\nabla(u-v)}_{L^2(K)}^2
&= \sum_{K\in\T_h}\inf_{v_K\in S(K)}\norm{\nabla(u-v_K)}_{L^2(K)}^2
\\&=\sum_{K\in\T_h}\inf_{w\in\hat S}\norm{\nabla(\hat u_K-w)}_{L^2(\hat K)}^2,
\end{aligned}
\end{equation}
where we have made the change of variable $w=\hat v_K$ in the last
step.

The proof proceeds by induction on $r$, the case $r=0$ being trivial.
For $r>0$,  apply \eqref{eq:sum1} with $u$ an arbitrary homogeneous
polynomial of degree $r$. Substituting \eqref{eq:hom} in
\eqref{eq:sum1}, and invoking the inductive hypothesis that $\P_0(\hat K)+\hat S\supseteq\P_{r-1}(\hat K)$, we get that
\begin{equation*}
\inf_{v\in S_h}\norm{\nabla_h(u-v)}_{L^2(\Omega)}^2
=h^{2r}\sum_{K\in\T_h}
\inf_{w\in\hat S}\norm{\nabla(u-w)}_{L^2(\hat K)}^2
=h^{2r-2}\inf_{w\in\hat S}\norm{\nabla(u-w)}_{L^2(\hat K)}^2.
\end{equation*}
Since we assume that this quantity is $o(h^{2r-2})$, the last infimum
must be $0$, so $u$ differs from an element $\hat S$ by a constant.
Thus $\P_0(\hat K)+\hat S$ contains all homogeneous
polynomials of degree $r$ and all polynomials of degree less than $r$
(by induction), so it indeed contains all polynomials of degree at most
$r$.
\end{proof}

\begin{remarks}
1.~ If $\hat S$ contains $\P_0(\hat K)$, which is usually the case,
then the third condition of Theorem~\ref{thm:poly1} reduces to that
of Theorem~\ref{thm:poly}.

2.~A similar result holds for higher derivatives (replace $\nabla_h$ by
$\nabla_h^m$ in the first two conditions, and
$\P_0(\hat K)$ by $\P_{m-1}(\hat K)$ in the third).
\end{remarks}

\section{Approximation theory of quadrilateral elements}
In this, the main section of the paper, we consider the approximation
properties of finite element spaces defined with respect to
quadrilateral meshes using bilinear mappings starting from a given
finite dimensional space of polynomials $\hat V$ on the unit square
$\hat K=[0,1]\x[0,1]$.  For simplicity we assume that $\hat
V\supseteq\P_0(\hat K)$. For example $\hat V$ might be the space
$\P_r(\hat K)$ of polynomials of degree at most $r$, or the space
$\Q_r(\hat K)$ of polynomials of degree at most $r$ in each variable
separately, or the serendipity space $\S_r(\hat K)$ spanned by
$\P_r(\hat K)$ together with the monomials $\hat x_1^r\hat x_2$ and
$\hat x_1\hat x_2^r$. Let $F$ be a bilinear isomorphism of $\hat K$
onto a convex quadrilateral $K=F(\hat K)$. Then for $u\in L^2(K)$ we
define $\hat u_F\in L^2(\hat K)$ by $\hat u_F(\hat x)=u(F\hat x)$, and
set
\begin{equation*}
V_F(K)=\{\, u:K\to\R\,|\, \hat u_F\in \hat V\,\}.
\end{equation*}
(Note that the definition of this space depends on the particular
choice of the bilinear isomorphism $F$ of $\hat K$ onto $K$, but whenever the space $\hat V$
is invariant under the symmetries of the square, which is usually the
case in practice, this will not be so.) We also note that the functions
in $V_F(K)$ need not be polynomials if $F$ is not affine, i.e., if $K$
is not a parallelogram.

Given a quadrilateral mesh $\T$ of some domain, $\Omega$, we can then
construct the space of functions $V^{\T}$ consisting of functions on the domain
which when restricted to a quadrilateral $K\in\T$ belong to
$V_{F_K}(K)$ where $F_K$ is a bilinear isomorphism of $\hat K$ onto $K$. 
(Again, if $\hat V$ is not invariant under the symmetries of the square,
the space $V^{\T}$ will depend on the specific choice of the maps $F_K$.)

It follows from the results of the previous section that if we consider
the sequence of meshes of the unit square into
congruent subsquares of side length $h=1/n$, then each of the
approximation estimates
\begin{gather}\label{eq:est1}
\inf_{v\in V^{\T_h}}\norm{u-v}_{L^2(\Omega)} \leq Ch^{r+1}\seminorm{u}_{r+1}
\text{\ for all $u\in H^{r+1}(\Omega)$},
\\\label{eq:est2}
\inf_{v\in V^{\T_h}}\norm{\nabla_h(u-v)}_{L^2(\Omega)} \leq
 Ch^r\seminorm{u}_{r+1}
\text{\ for all $u\in H^{r+1}(\Omega)$}, 
\end{gather}
holds if and only $\P_r(\hat K)\subset \hat V$.  It is not hard to
extend these estimates to shape-regular sequences of parallelogram meshes
as well.
\emph{However, in this section we show that for these estimates to hold
for more general quadrilateral mesh sequences, a stronger condition
on $\hat V$ is required, namely that $\hat V\supseteq\Q_r(\hat K)$.}

The positive result, that when $\hat V\supseteq\Q_r(\hat K)$, then
the estimates \eqref{eq:est1} and \eqref{eq:est2} hold
for any shape regular sequence of quadrilateral meshes $\T_h$,
is known.  See, e.g., \cite{ciarlet-raviart}, \cite{ciarlet}, or
\cite[Section I.A.2]{girault-raviart}.
We wish to show the necessity of the condition $\hat V\supseteq\Q_r(\hat K)$.

As a first step we show that the condition $V_F(K)\supseteq\P_r(K)$ is
necessary and sufficient to have that $\hat V\supseteq\Q_r(\hat K)$
whenever $F$ is a bilinear isomorphism of $\hat K$ onto a convex
quadrilateral.  This is proven in the following two theorems.

\begin{thm}
Suppose that $\hat V\supseteq\Q_r(\hat K)$.  Let $F$ be any bilinear
isomorphism of $\hat K$ onto a convex quadrilateral.  Then
$V_F(K)\supseteq\P_r(K)$.
\end{thm}
\begin{proof}
The components of $F(\hat x,\hat y)$ are linear functions of $\hat x$
and $\hat y$, so if $p$ is a polynomial of degree at most $r$, then
$p(F(\hat x,\hat y))$ is of degree at most $r$ in $\hat x$ and $\hat y$,
i.e., $p\circ F\in\Q_r(\hat K)\subset \hat V$. Therefore $p\in V_F(K)$.
\end{proof}

The reverse implication holds even under the weaker assumption
that $V_F(K)$ contains $\P_r(K)$ just for the two specific 
bilinear isomorphism
\begin{equation*}
\tilde F(\hat x,\hat y)=(\hat x,\hat y(\hat x+1)),\quad
\bar F(\hat x,\hat y)=(\hat y,\hat x(\hat y+1)),
\end{equation*}
both of which map $\hat K$ isomorphically onto the quadrilateral
$K'$ with vertices $(0,0)$, $(1,0)$, $(0,1)$, and $(1,2)$.
This fact is established below.

\begin{thm}\label{thm:nopoly}
Let $\hat V$ be a vectorspace of functions on $\hat K$. Suppose that
$\Q_r(\hat K)\nsubseteq\hat V$.  Then either $V_{\tilde F}(K')\nsubseteq\P_r(K')$ or
$V_{\bar F}(K')\nsubseteq\P_r(K')$. 
\end{thm}

\begin{remark}
If the space $\hat V$ is invariant under the symmetries of the square, then
$V_{\tilde F}(K') = V_{\bar F}(K')$ so neither contains $\P_r(K')$.
\end{remark}

\begin{proof}
Assume to the contrary that $V_{\tilde F}(K')\supseteq\P_r(K')$ and
$V_{\bar F}(K')\supseteq\P_r(K')$. We prove that $\hat V\supseteq\Q_r(\hat K)$ by induction on
$r$.  The case $r=0$ being true by assumption, we consider $r>0$ and show
that the monomials $\hat x^r\hat y^s$ and $\hat x^s\hat y^r$ belong to
$\hat V$ for $s=0,1,\ldots,r$.  From the identity
\begin{equation}\label{eq:id}
\hat x^r\hat y^s = \hat x^{r-s}[\hat y(\hat x+1)]^s
-\sum_{t=1}^s\binom st \hat x^{r-t}\hat y^s
= \tilde F_1(\hat x,\hat y)^{r-s}\tilde F_2(\hat x,\hat y)^s
-\sum_{t=1}^s\binom st \hat x^{r-t}\hat y^s,
\end{equation}
we see that for $0\le s<r$, the monomial $\hat x^r\hat y^s$ is the sum
of a polynomial which clearly belongs to $\hat V$ (since $
\tilde F_1(\hat x,\hat y)^{r-s}\tilde F_2(\hat x,\hat y)^s = x^{r-s} y^s
\in \P_r(K') \subset V_{\tilde F}(K')$) and a polynomial
in $Q_{r-1}(\hat K)$, which belongs to $\hat V$ by induction.
Thus each of the monomials $\hat x^r\hat y^s$ with $0\le s<r$ belongs
to $\hat V$, and, using $\bar F$, we similarly see that all the monomials
$\hat x^s\hat y^r$, $0\le s<r$ belong to $\hat V$.  Finally, from
\eqref{eq:id} with $s=r$, we see that $\hat x^r\hat y^r$ is a linear
combination of an element of $\hat V$ and monomials
$\hat x^s\hat y^r$ with $s<r$, so it too belongs to $\hat V$.
\end{proof}

\begin{figure}[ht]
\centerline{\includegraphics[height=1.5in]{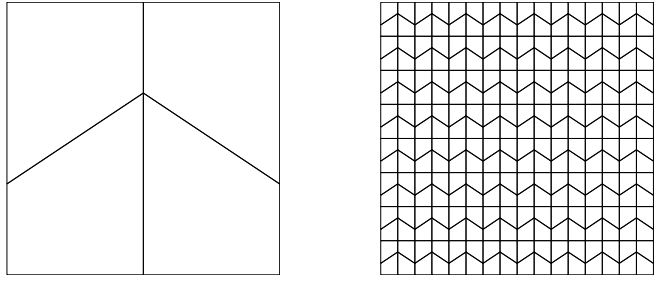}}
\caption{\label{fg:macro} a.~A partition of the square into
four trapezoids.  b.~A mesh composed of translated dilates of this
partition.}
\end{figure}

We now combine this result with the those of the previous section to
show the necessity of the condition  $\hat V\supseteq\Q_r(\hat K)$
for optimal order approximation.
Let $\hat V$ be some fixed finite dimensional subspace of $L^2(\hat K)$
which does \emph{not} include $Q_r(\hat K)$.
Consider the specific division of the unit square $\hat K$ into
four quadrilaterals shown on the left in Figure~\ref{fg:macro}.  For definiteness
we place the vertices of the quadrilaterals at $(0,1/3)$, $(1/2,2/3)$
and $(1,1/3)$ and the midpoints of the horizontal edges and the corners of $\hat K$.

\medskip\nobreak
The four quadrilaterals are mutually congruent and affinely
related to the specific quadrilateral $K'$ defined above.  Therefore, by
Theorem~\ref{thm:nopoly},  we can define for each of the four
quadrilaterals $K''$ shown in Figure~\ref{fg:macro} an isomorphism $F''$ from
the unit square so that $V_{F''}(K'')\nsupseteq\P_r(K'')$. If we let $\hat S$ be
the subspace of $L^2(\hat K)$ consisting of functions which restrict
to elements of $V_{F''}(K'')$ on each of the four quadrilaterals $K''$,
then certainly $\hat S$ does not contain $\P_r(\hat K)$, since
even its restriction to any one of the quadrilaterals $K''$ does not
contain $P_r(K'')$.

Next, for $n=1,2,\ldots$ consider the mesh $\T'_h$ of the unit square
$\Omega$ shown in Figure~\ref{fg:macro}b, obtained by first dividing
it into a uniform $n\x n$ mesh of subsquares, $n=1/h$, and then
dividing each subsquare as in Figure~\ref{fg:macro}a.
Then the space of functions $u$ on $\Omega$ whose restrictions on each
subsquare $K\in\T_h$ satisfy $\hat u_K(\hat x)=u(x_K+h\hat x)$ with
$\hat u_K\in\hat S$ is precisely the same as the space $V(\T'_h)$
constructed from the initial space $\hat V$ and the mesh $\T'_h$. In
view of Theorems~\ref{thm:poly} and \ref{thm:poly1} and the fact that
$\hat S\nsupseteq\P_r(\hat K)$, the estimates \eqref{eq:est1} and
\eqref{eq:est2} do not hold.  In fact, neither of the estimates
\begin{equation*}
\inf_{v\in V(\T_h)}\norm{u-v}_{L^2(\Omega)} = o(h^r),
\end{equation*}
nor
\begin{equation*}
\inf_{v\in V(\T_h)}\norm{\nabla(u-v)}_{L^2(\Omega)} = o(h^{r-1}),
\end{equation*}
holds, even for $u\in\P_r(\Omega)$.

While the condition $\hat V\supseteq \Q_r(\hat K)$ is necessary for
$O(h^{r+1})$ on general quadrilateral meshes, the conditions $\hat
V\supseteq \P_r(\hat K)$ suffices for meshes of parallelograms. 
Naturally, the same is true for meshes whose elements are sufficiently
close to parallelograms.  We conclude this section with a precise
statement of this result and a sketch of the proof. If $\hat
V\supseteq \P_r(\hat K)$ and $K=F(\hat K)$ with $F\in\Bil(\hat K)$,
then by standard arguments, as in \cite{ciarlet}, we get
\begin{equation*}
\|v-\pi v\|_{L^2(K)} \le C \|J_F\|_{L^\infty(\hat K)}^{1/2}|v\circ
F|_{H^{r+1}(\hat K)},
\end{equation*}
where $J_F$ is the Jacobian determinant of $F$.
Now, using the formula for the derivative of a composition 
(as in, e.g., \cite[p.~222]{Federer}), and the
fact that $F$ is quadratic, and so its third and higher derivatives
vanish, we get
that
\begin{equation*}
|v\circ F|_{H^{r+1}(\hat K)}\le C\|J_{F^{-1}}\|_{L^\infty(K)}^{1/2}
\|v\|_{H^{r+1}(K)}\sum_{i=0}^{\lfloor (r+1)/2\rfloor} |F|_{W^1_\infty(\hat
K)}^{r+1-2i}|F|_{W^2_\infty(\hat K)}^i.
\end{equation*}
Now,
\begin{equation*}
\|J_F\|_{L^\infty(\hat K)}\le Ch_K^2, \quad
\|J_{F^{-1}}\|_{L^\infty(\hat K)}\le Ch_K^{-2},
\quad
|F|_{W^1_\infty(\hat K)}\le C h_K,
\end{equation*}
where $h_K$ is the diameter of $K$ and $C$ depends only on the
shape-regularity of $K$.  We thus get
\begin{equation*}
\|v-\pi v\|_{L^2(K)} \le C \|v\|_{H^{r+1}(K)}\sum_i h_K^{r+1-2i}|F|_{W^2_\infty(\hat K)}^i.
\end{equation*}
It follows that if $|F|_{W^2_\infty(\hat K)}=O(h_K^2)$, we get the
desired estimate
\begin{equation*}
\|v-\pi v\|_{L^2(K)} \le C h_K^{r+1}\|v\|_{H^{r+1}(K)}.
\end{equation*}

Following \cite{rannacher-turek}, we measure the deviation of a
quadrilateral from a parallelogram, by the quantity
$\sigma_K:=\max(|\pi-\theta_1|,|\pi-\theta_2|)$, where $\theta_1$ is the
angle between the outward normals of two opposite sides of $K$ and
$\theta_2$ is the angle between the outward normals of the other two
sides.  Thus $0\le \sigma_K<\pi$, with $\sigma_K=0$ if and only if $K$
is a parallelogram.  As pointed out in \cite{rannacher-turek},
$|F|_{W^2_\infty(\hat K)}\le Ch_K(h_K+\sigma_K)$.  This motivates
the definition that a family of quadrilateral meshes is asymptotically
parallelogram if $\sigma_K=O(h_K)$, i.e., if $\sigma_K/h_K$ is uniformly
bounded for all the elements in all the meshes.  From the foregoing
considerations, if the reference space contains $\P_r(\hat K)$ we
obtain $O(h^{r+1})$ convergence for asymptotically parallelogram, shape regular
meshes.

As a final note, we remark that any polygon can be meshed by an
asymptotically parallelogram, shape regular family of meshes with mesh size
tending to zero.  Indeed, if we begin with any mesh of convex
quadrilaterals, and refine it by dividing each quadrilateral in four by
connecting the midpoints of the opposite edges, and continue in this
fashion, as in the last row of Figure~\ref{fg:meshes}, the resulting
mesh is asymptotically parallelogram and shape regular.

\section{Numerical results}
\begin{figure}[ht]
\centerline{\includegraphics[width=4in]{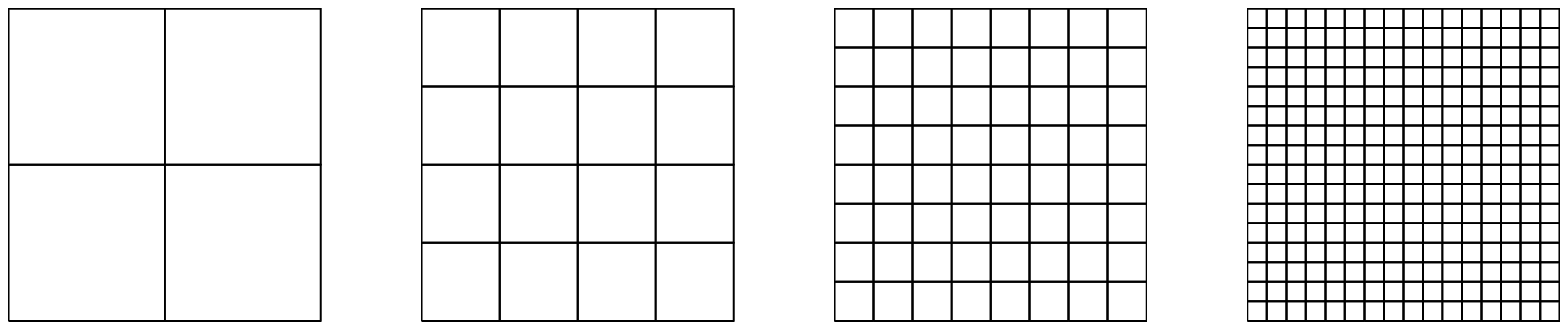}}
\centerline{\includegraphics[width=4in]{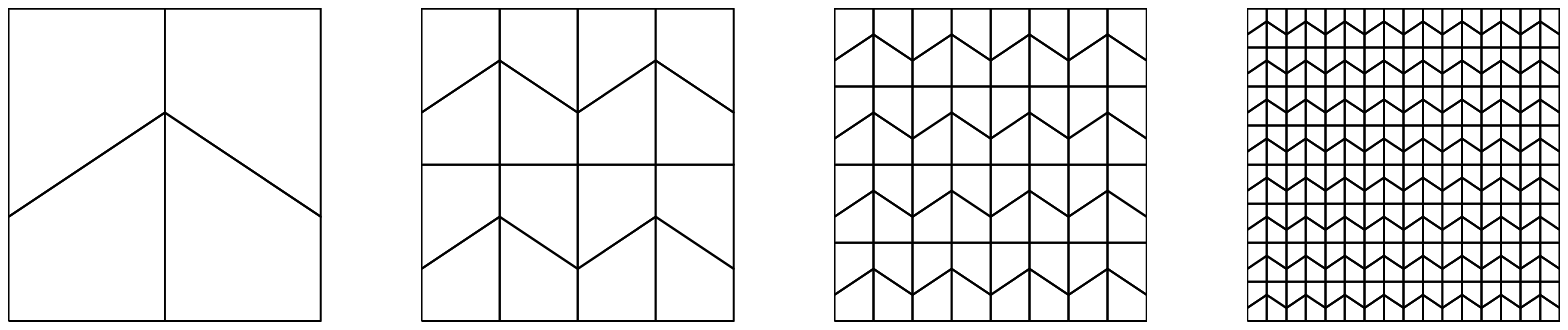}}
\centerline{\includegraphics[width=4in]{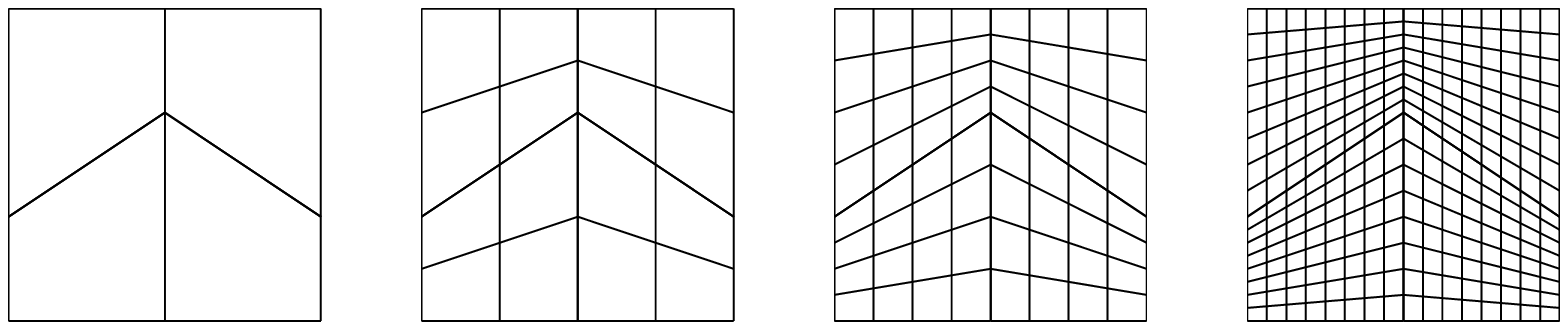}}
\caption{\label{fg:meshes} \footnotesize Three sequences of meshes of the unit square:
square, trapezoidal, and asymptotically parallelogram.
Each is shown for $n=2$ ,$4$, $8$, and $16$.}
\end{figure}

In this section we report on results from a numerical study of the
behavior of piecewise continuous mapped biquadratic and serendipity
finite elements on quadrilateral meshes (i.e., the finite element spaces
are constructed starting from the spaces $Q_2(\hat K)$ and $S_2(\hat K)$
on the reference square, and then imposing continuity).
We present the results of two
test problems.  In both we solve the Dirichlet problem for Poisson's
equation
\begin{equation}\label{eq:prob}
-\Delta u = f \text{ in $\Omega$},\quad
u = g \text{ on $\partial\Omega$},
\end{equation}
where the domain $\Omega$ is the unit square.  In the first problem,
$f$ and $g$ are taken so that the exact solution is the quartic
polynomial
\begin{equation*}
u(x,y)=x^3+5 y^2 - 10y^3 + y^4.
\end{equation*}
Table~\ref{tb:t1} shows results for both types of elements using meshes
from each of the first two mesh sequences shown in
Figure~\ref{fg:meshes}. The first sequence of meshes consists of
uniform square subdivisions of the domain into $n\x n$ subsquares,
$n=2,4,8,\ldots$.  Meshes in the second sequence are partitions of the
domain into $n\x n$ congruent trapezoids, all similar to the trapezoid
with vertices $(0,0)$, $(1/2,0)$, $(1/2,2/3)$, and $(0,1/3)$. In
Table~\ref{tb:t1} we report the errors in $L^2$ for the finite element
solution and its gradient both in absolute terms and as a percentage of
the $L^2$ norm of the exact solution and its gradient, and we also
report the apparent rate of convergence based on consecutive meshes in
a sequence. For this test problem, the rates of convergence are very
clear: for either mesh sequence the mapped biquadratic elements
converge with the expected order $3$ for the solution and $2$ for its
gradient.  The same is true for the serendipity elements on the square
meshes, but, as predicted by the theory given above, for the
trapezoidal mesh sequence the order of convergence for the serendipity
elements is reduced by $1$ both for the solution and its gradient. 

%Table 1
\begin{table}[tbq]
\footnotesize
\caption{ Errors and rates of convergence for the test problem with polynomial solution.}
\label{tb:t1}
\begin{center}
\begin{tabular}{r|rrr|rrr|rrr|rrr}
\multicolumn{13}{c}{\emph{Mapped biquadratic  elements}} \\[.5ex]
\multicolumn{1}{c}{} & \multicolumn{6}{c}{square meshes} & \multicolumn{6}{c}{trapezoidal meshes} \\
\hline
\rule{0pt}{15pt}& \multicolumn{3}{c|}{$\|u-u_h\|_{L^2}$} & 
\multicolumn{3}{c|}{$\|\nabla(u-u_h)\|_{L^2}$} &
\multicolumn{3}{c|}{$\|u-u_h\|_{L^2}$} & 
\multicolumn{3}{c}{$\|\nabla(u-u_h)\|_{L^2}$} \\[.5ex]
\multicolumn{1}{c|}{$n$} &
\multicolumn{1}{c}{err.} & 
\multicolumn{1}{c}{\%} & 
\multicolumn{1}{c|}{rate} & 
\multicolumn{1}{c}{err.} & 
\multicolumn{1}{c}{\%} & 
\multicolumn{1}{c|}{rate} & 
\multicolumn{1}{c}{err.} & 
\multicolumn{1}{c}{\%} & 
\multicolumn{1}{c|}{rate} & 
\multicolumn{1}{c}{err.} & 
\multicolumn{1}{c}{\%} & 
\multicolumn{1}{c}{rate} \\[.5ex]
\hline
\rule{0pt}{15pt}
  2 & 3.5e$-$02 &   2.877 &     & 4.5e$-$01 &  37.253 &     & 4.8e$-$02 &   3.951 &     & 5.9e$-$01 &  48.576 &	    \\
  4 & 4.4e$-$03 &   0.360 & 3.0 & 1.1e$-$01 &   9.333 & 2.0 & 5.8e$-$03 &   0.475 & 3.1 & 1.5e$-$01 &  12.082 & 2.0 \\
  8 & 5.5e$-$04 &   0.045 & 3.0 & 2.8e$-$02 &   2.329 & 2.0 & 7.1e$-$04 &   0.058 & 3.0 & 3.7e$-$02 &   3.017 & 2.0 \\
 16 & 6.9e$-$05 &   0.006 & 3.0 & 7.1e$-$03 &   0.583 & 2.0 & 8.7e$-$05 &   0.007 & 3.0 & 9.2e$-$03 &   0.753 & 2.0 \\
 32 & 8.6e$-$06 &   0.001 & 3.0 & 1.8e$-$03 &   0.146 & 2.0 & 1.1e$-$05 &   0.001 & 3.0 & 2.3e$-$03 &   0.188 & 2.0 \\
 64 & 1.1e$-$06 &   0.000 & 3.0 & 4.4e$-$04 &   0.036 & 2.0 & 1.3e$-$06 &   0.000 & 3.0 & 5.7e$-$04 &   0.047 & 2.0 \\
%128 & 1.4e$-$07 &   0.000 & 2.9 & 1.1e$-$04 &   0.009 & 2.0 & 1.9e$-$07 &   0.000 & 2.9 & 1.4e$-$04 &   0.012 & 2.0 \\[.5ex]
\hline
\multicolumn{13}{c}{\rule{0pt}{20pt}\emph{Serendipity elements}}\\[.5ex]
\multicolumn{1}{c}{} & \multicolumn{6}{c}{square meshes} & \multicolumn{6}{c}{trapezoidal meshes} \\
\hline
\rule{0pt}{15pt}& \multicolumn{3}{c|}{$\|u-u_h\|_{L^2}$} & 
\multicolumn{3}{c|}{$\|\nabla(u-u_h)\|_{L^2}$} &
\multicolumn{3}{c|}{$\|u-u_h\|_{L^2}$} & 
\multicolumn{3}{c}{$\|\nabla(u-u_h)\|_{L^2}$} \\[.5ex]
\multicolumn{1}{c|}{$n$} &
\multicolumn{1}{c}{err.} & 
\multicolumn{1}{c}{\%} & 
\multicolumn{1}{c|}{rate} & 
\multicolumn{1}{c}{err.} & 
\multicolumn{1}{c}{\%} & 
\multicolumn{1}{c|}{rate} & 
\multicolumn{1}{c}{err.} & 
\multicolumn{1}{c}{\%} & 
\multicolumn{1}{c|}{rate} & 
\multicolumn{1}{c}{err.} & 
\multicolumn{1}{c}{\%} & 
\multicolumn{1}{c}{rate} \\[.5ex]
\hline
\rule{0pt}{15pt}
  2 & 3.5e$-$02 &   2.877 &     & 4.5e$-$01 &  37.252 &     & 5.0e$-$02 &   4.066 &     & 6.2e$-$01 &  51.214 &	    \\
  4 & 4.4e$-$03 &   0.360 & 3.0 & 1.1e$-$01 &   9.333 & 2.0 & 6.7e$-$03 &   0.548 & 2.9 & 1.8e$-$01 &  14.718 & 1.8 \\
  8 & 5.5e$-$04 &   0.045 & 3.0 & 2.8e$-$02 &   2.329 & 2.0 & 9.7e$-$04 &   0.080 & 2.8 & 5.9e$-$02 &   4.836 & 1.6 \\
 16 & 6.9e$-$05 &   0.006 & 3.0 & 7.1e$-$03 &   0.583 & 2.0 & 1.6e$-$04 &   0.013 & 2.6 & 2.3e$-$02 &   1.890 & 1.4 \\
 32 & 8.6e$-$06 &   0.001 & 3.0 & 1.8e$-$03 &   0.146 & 2.0 & 3.3e$-$05 &   0.003 & 2.3 & 1.0e$-$02 &   0.842 & 1.2 \\
 64 & 1.1e$-$06 &   0.000 & 3.0 & 4.4e$-$04 &   0.036 & 2.0 & 7.4e$-$06 &   0.001 & 2.1 & 4.9e$-$03 &   0.401 & 1.1 \\
%128 & 1.4e$-$07 &   0.000 & 2.9 & 1.1e$-$04 &   0.009 & 2.0 & 1.8e$-$06 &   0.000 & 2.1 & 2.4e$-$03 &   0.197 & 1.0 \\[.5ex]
\hline
\end{tabular}
\end{center}
\end{table}

As a second test example we again solved the Dirichlet problem
\eqref{eq:prob}, but this time choosing the data so that the solution
is the sharply peaked function
\begin{equation*}
u(x,y)=\exp\bigl(-100[(x-1/4)^2+(y-1/3)^2]\bigr).
\end{equation*}
As seen in Table~\ref{tb:t2}, in this case the loss of convergence
order for the serendipity elements on the trapezoidal mesh is not
nearly as clear.  Some loss is evident, but apparently very fine meshes
(and very high precision computation) would be required to see the
final asymptotic orders.

%Table 2
\begin{table}[tbp]
\footnotesize
\caption{ Errors and rates of convergence for the test problem with
exponential solution.}
\label{tb:t2}
\begin{center}
\begin{tabular}{r|rrr|rrr|rrr|rrr}
\multicolumn{13}{c}{\emph{Mapped biquadratic  elements}} \\[.5ex]
\multicolumn{1}{c}{} & \multicolumn{6}{c}{square meshes} & \multicolumn{6}{c}{trapezoidal meshes} \\
\hline
\rule{0pt}{15pt}& \multicolumn{3}{c|}{$\|u-u_h\|_{L^2}$} & 
\multicolumn{3}{c|}{$\|\nabla(u-u_h)\|_{L^2}$} &
\multicolumn{3}{c|}{$\|u-u_h\|_{L^2}$} & 
\multicolumn{3}{c}{$\|\nabla(u-u_h)\|_{L^2}$} \\[.5ex]
\multicolumn{1}{c|}{$n$} &
\multicolumn{1}{c}{err.} & 
\multicolumn{1}{c}{\%} & 
\multicolumn{1}{c|}{rate} & 
\multicolumn{1}{c}{err.} & 
\multicolumn{1}{c}{\%} & 
\multicolumn{1}{c|}{rate} & 
\multicolumn{1}{c}{err.} & 
\multicolumn{1}{c}{\%} & 
\multicolumn{1}{c|}{rate} & 
\multicolumn{1}{c}{err.} & 
\multicolumn{1}{c}{\%} & 
\multicolumn{1}{c}{rate} \\[.5ex]
\hline
\rule{0pt}{15pt}
  2 & 2.8e$-$01 & \llap{2}24.000 &     & 3.0e$+$00 & \llap{1}69.630 &    
  & 2.6e$-$01 & \llap{2}04.800 &     & 2.8e$+$00 & \llap{1}59.208 &	    \\
  4 & 1.2e$-$01 &  93.600 & 1.3 & 1.5e$+$00 &  87.322 & 1.0 & 2.1e$-$01 & 169.600 & 0.3 & 1.8e$+$00 &  99.305 & 0.7 \\
  8 & 1.7e$-$02 &  13.520 & 2.8 & 4.6e$-$01 &  25.809 & 1.8 & 2.3e$-$02 &  18.160 & 3.2 & 5.9e$-$01 &  33.185 & 1.6 \\
 16 & 1.1e$-$03 &   0.920 & 3.9 & 1.0e$-$01 &   5.860 & 2.1 & 1.3e$-$03 &   1.048 & 4.1 & 1.2e$-$01 &   6.819 & 2.3 \\
 32 & 1.3e$-$04 &   0.101 & 3.2 & 2.5e$-$02 &   1.424 & 2.0 & 1.5e$-$04 &   0.124 & 3.1 & 3.2e$-$02 &   1.794 & 1.9 \\
 64 & 1.5e$-$05 &   0.012 & 3.1 & 6.3e$-$03 &   0.354 & 2.0 & 1.9e$-$05 &   0.015 & 3.0 & 7.9e$-$03 &   0.448 & 2.0 \\
128 & 1.9e$-$06 &   0.002 & 3.0 & 1.6e$-$03 &   0.088 & 2.0 & 2.4e$-$06 &   0.002 & 3.0 & 2.0e$-$03 &   0.112 & 2.0 \\[.5ex]
\hline
\multicolumn{13}{c}{\rule{0pt}{20pt}\emph{Serendipity elements}}\\[.5ex]
\multicolumn{1}{c}{} & \multicolumn{6}{c}{square meshes} & \multicolumn{6}{c}{trapezoidal meshes} \\
\hline
\rule{0pt}{15pt}& \multicolumn{3}{c|}{$\|u-u_h\|_{L^2}$} & 
\multicolumn{3}{c|}{$\|\nabla(u-u_h)\|_{L^2}$} &
\multicolumn{3}{c|}{$\|u-u_h\|_{L^2}$} & 
\multicolumn{3}{c}{$\|\nabla(u-u_h)\|_{L^2}$} \\[.5ex]
\multicolumn{1}{c|}{$n$} &
\multicolumn{1}{c}{err.} & 
\multicolumn{1}{c}{\%} & 
\multicolumn{1}{c|}{rate} & 
\multicolumn{1}{c}{err.} & 
\multicolumn{1}{c}{\%} & 
\multicolumn{1}{c|}{rate} &
\multicolumn{1}{c}{err.} & 
\multicolumn{1}{c}{\%} & 
\multicolumn{1}{c|}{rate} & 
\multicolumn{1}{c}{err.} & 
\multicolumn{1}{c}{\%} & 
\multicolumn{1}{c}{rate} \\[.5ex]
\hline
\rule{0pt}{15pt}
  2 & 2.0e$-$01 & \llap{1}59.200 &     & 2.4e$+$00 & \llap{1}33.372 &     & 2.1e$-$01 & \llap{1}69.600 &     & 2.3e$+$00 & \llap{1}30.340 &	    \\
  4 & 1.2e$-$01 &  92.000 & 0.8 & 1.4e$+$00 &  80.531 & 0.7 & 2.1e$-$01 & 168.000 & 0.0 & 1.7e$+$00 &  93.819 & 0.5 \\
  8 & 1.7e$-$02 &  13.520 & 2.8 & 4.6e$-$01 &  26.293 & 1.6 & 2.4e$-$02 &  18.880 & 3.2 & 6.1e$-$01 &  34.564 & 1.4 \\
 16 & 1.1e$-$03 &   0.920 & 3.9 & 1.1e$-$01 &   5.948 & 2.1 & 1.5e$-$03 &   1.208 & 4.0 & 1.4e$-$01 &   7.737 & 2.2 \\
 32 & 1.3e$-$04 &   0.101 & 3.2 & 2.5e$-$02 &   1.432 & 2.1 & 2.0e$-$04 &   0.162 & 2.9 & 3.8e$-$02 &   2.156 & 1.8 \\
 64 & 1.5e$-$05 &   0.012 & 3.1 & 6.3e$-$03 &   0.354 & 2.0 & 2.7e$-$05 &   0.022 & 2.9 & 1.1e$-$02 &   0.597 & 1.9 \\
128 & 1.9e$-$06 &   0.002 & 3.0 & 1.6e$-$03 &   0.088 & 2.0 & 3.7e$-$06 &   0.003 & 2.9 & 3.4e$-$03 &   0.191 & 1.6 \\[.5ex]
\hline
\end{tabular}
\end{center}
\end{table}

Finally we return to the first test problem, and consider the behavior
of the serendipity elements on the third mesh sequence shown in
Figure~\ref{fg:meshes}.  This mesh sequence begins with the same mesh
of four quadrilaterals as in previous case, and continues with
systematic refinement as described at the end of the last section, and
so is asymptotically parallelogram.  Therefore, as explained there, the
rate of convergence for serendipity elements is the
same as for affine meshes.   This is clearly illustrated in
Table~\ref{tb:t3}.

%Table 3
\begin{table}[tbp]
\footnotesize
\caption{ Errors and rates of convergence for the test problem with
polynomial solution using serendipity elements on asympotically affine
meshes.}
\label{tb:t3}
\begin{center}
\begin{tabular}{r|rrr|rrr}
\hline
\rule{0pt}{15pt}& \multicolumn{3}{c|}{$\|u-u_h\|_{L^2}$} & 
\multicolumn{3}{c}{$\|\nabla(u-u_h)\|_{L^2}$} \\[.5ex]
\multicolumn{1}{c|}{$n$} &
\multicolumn{1}{c}{err.} & 
\multicolumn{1}{c}{\%} & 
\multicolumn{1}{c|}{rate} & 
\multicolumn{1}{c}{err.} & 
\multicolumn{1}{c}{\%} & 
\multicolumn{1}{c}{rate} \\[.5ex]
\hline
\rule{0pt}{15pt}
  2 & 5.0e$-$02 &   4.066 &     & 6.2e$-$01 &  51.214 &     \\
  4 & 6.2e$-$03 &   0.510 & 3.0 & 1.5e$-$01 &  12.109 & 2.1 \\
  8 & 7.6e$-$04 &   0.062 & 3.0 & 3.6e$-$02 &   2.948 & 2.0 \\
 16 & 9.4e$-$05 &   0.008 & 3.0 & 9.0e$-$03 &   0.735 & 2.0 \\
 32 & 1.2e$-$05 &   0.001 & 3.0 & 2.2e$-$03 &   0.183 & 2.0 \\
 64 & 1.5e$-$06 &   0.000 & 3.0 & 5.6e$-$04 &   0.046 & 2.0 \\
128 & 1.9e$-$07 &   0.000 & 3.0 & 1.4e$-$04 &   0.012 & 2.0 \\
\hline
\end{tabular}
\end{center}
\end{table}

While the asymptotic rates predicted by the theory are confirmed in
these examples, it is worth noting that in absolute terms the effect of
the degraded convergence rate is not very pronounced.  For the first
example, on a moderately fine mesh of $16\x 16$ trapezoids, the
solution error with serendipity elements exceeds that of mapped
biquadratic elements by a factor of about 2, and the gradient error by
a factor of 2.5.  Even on the finest mesh shown, with $64\x64$
elements, the factors are only about 5.5 and 8.5, respectively. Of
course, if we were to compute on finer and finer meshes with
sufficiently high precision, these factors would tend to infinity. 
Indeed, on any quadrilateral mesh which contains a non-parallelogram
element, the analogous factors can be made as large as desired by
choosing a problem in which the exact solution is sufficiently close
to---or even equal to---a quadratic function, which the mapped
biquadratic elements capture exactly, while the serendipity elements do
not (such a quadratic function always exists). However, it is not
unusual that the serendipity elements perform almost as well as the
mapped biquadratic elements for reasonable, and even for quite small,
levels of error.  This, together with their optimal convergence on
asymptotically parallelogram meshes, provides an explanation of why the
lower rates of convergence have not been widely noted.

\bibliographystyle{amsplain}
\bibliography{quadapprox}

\end{document}